\documentclass[12pt]{amsart2000}
\usepackage{amsmath2000,amssymb,amsfonts,amsthm2000,epsfig,xspace}
\usepackage[round]{natbib}
\newtheorem{theorem}{Theorem}
\newtheorem{lemma}[theorem]{Lemma}

\newcommand{\tref}[1]{Theorem~\ref{thm:#1}}
\newcommand{\lref}[1]{Lemma~\ref{lem:#1}}

\newcommand{\E}{{\mathbb E}}
\newcommand{\eps}{\varepsilon}

\newcommand{\dsc}{Diaconis and Saloff-Coste\xspace}
\newcommand{\gap}{\gamma}

\newcommand{\dps}{\Delta\Psi}

\newcommand{\px}{\Psi_{\max}}

\newcommand{\Var}{\operatorname{Var}}

\newcommand{\TV}{\operatorname{TV}}

\begin{document}

\title{Mixing Time of the Rudvalis Shuffle}
\author[David B. Wilson]{David Bruce Wilson}
\address{Microsoft Research\\One Microsoft Way\\Redmond, WA 98052\\U.S.A.}
\email{dbwilson@microsoft.com}
\urladdr{http://dbwilson.com}
\date{}
\begin{abstract}
We extend a technique for lower-bounding the mixing time of
card-shuffling Markov chains, and use it to bound the mixing time of
the Rudvalis Markov chain, as well as two variants considered by
Diaconis and Saloff-Coste.  We show that in each case $\Theta(n^3\log
n)$ shuffles are required for the permutation to randomize, which
matches (up to constants) previously known upper bounds.  In contrast,
for the two variants, the mixing time of an individual card is only
$\Theta(n^2)$ shuffles.
\end{abstract}
\maketitle
\vspace*{-12pt}

\section{Introduction}

In earlier work \citep{wilson:mixing-times} we derived upper and lower
bounds on the mixing time of a variety of Markov chains, including
Markov chains on lozenge tilings, card shuffling, and exclusion
processes.  The mixing time of a Markov chain is the time it takes to
approach its stationary distribution, which is often measured in total
variation distance (defined below).  In this article we focus on
the method for lower bounding the mixing time, and extend its
applicability to the Rudvalis card shuffling Markov chain (defined
below) and related shuffles.

Let $P^{*t}_x$ denote the distribution of the Markov chain started in
state $x$ after it is run for $t$ steps, and let $\mu$ denote the
stationary distribution of the Markov chain.  The total variation
distance between distributions $P^{*t}$ and $\mu$ is defined by
$$ \left\|P^{*t}_x-\mu\right\|_{\TV} = \max_A\left|P^{*t}_x(A)-\mu(A)\right| = \frac{1}{2} \sum_y \left|P^{*t}_x(y)-\mu(y)\right| = \frac{1}{2} \left\|P^{*t}-\mu\right\|_1,$$
and the mixing time is the time it takes for
$\max_x \left\|P^{*t}_x-\mu\right\|_{\TV}$ to become small,
say smaller than $\eps$.
See \citep{aldous-fill:book,diaconis:monograph,diaconis:cutoff} for further background.

Arunas Rudvalis proposed the following shuffle: with probability $1/2$
move the top card to the bottom of the deck, and with probability
$1/2$ move it to the second position from the bottom.  \cite{hildebrand:thesis}
showed that the Rudvalis shuffle mixes in $O(n^3\log n)$ time.  
\cite{diaconis-saloff-coste:survey} studied a variation,
the shift-or-swap shuffle, which at each step either moves the top
card to the bottom of the deck or exchanges the top two cards, each
move with probability $1/2$.  \cite{diaconis-saloff-coste:comparison}
also studied a symmetrized version of the Rudvalis shuffle, which at
each step does one of four moves each with probability $1/4$: move top
card to bottom, move bottom card to top, exchange top two cards, or do
nothing.  In each case a $O(n^3\log n)$ upper bound on the mixing time
is known, but order $n^3\log n$ lower bounds were not known.

To lower bound the mixing time, one finds a set $A$ of states such
that $P^{*t}(A)$ is close to $1$ and $\mu(A)$ is close to $0$.  The
approach taken in \citep{wilson:mixing-times}
uses an eigenvector $\Phi$ of the Markov
chain.  If $X_t$ denotes the state of the Markov chain at time $t$,
then $\E[\Phi(X_{t+1})\mid X_t] = \lambda \Phi(X_t)$.  To obtain a good
lower bound, we need $\lambda<1$ but $\lambda\approx 1$.  Since
$\lambda<1$, in stationarity $\E[\Phi(X)]=0$, but since
$\lambda\approx 1$, it takes a long time before $\E[\Phi(X_t)]\approx
0$.  If furthermore the eigenvector is ``smooth'' in the sense that
$\E[|\Phi(X_{t+1})-\Phi(X_t)|^2\mid X_t]$ is never large, then we can
bound the variance $\Phi(X_t)$, showing that it is with high
probability confined to a small interval about its expected value.  Then
provided that $\E[\Phi(X_t)]$ is large enough, we can reliably
distinguish $\Phi(X_t)$ from $\Phi(X)$ in stationarity, which implies
that the Markov chain has not yet mixed by time $t$.
\cite{saloff-coste:lower-bound} gives an exposition of this and
related ideas.

For the Rudvalis card shuffling Markov chain and its variants,
there are a few difficulties when directly applying this
approach to lower bound the mixing time.  The eigenvectors that one
wants to use are complex-valued rather than real-valued, and
$\Phi(X_t)$ is no longer confined to a small interval around
$\E[\Phi(X_t)]$.  Instead what happens is that
$\Phi(X_t)$ is with high probability confined to a narrow annulus
centered at $0$.  $\E[\Phi(X_t)]$ becomes too small too quickly, and
$\Var[\Phi(X_t)]$ remains too large to be useful.

To lower bound the mixing time we want to in effect work with
$|\Phi(X_t)|$ and forget about $\arg\Phi(X_t)$.  To do this we start
by lifting the Markov chain to a larger state space, and let us denote
the state at time $t$ of the lifted chain by $(X_t,Y_t)$.  (The mixing
time of the lifted Markov chain will upper bound the mixing time of
the original chain, so at the outset it is not clear that we can lower
bound the mixing time of the original chain by considering its lifted
version.)  We find an eigenvector $\Psi$ on the lifted chain such that
for all $x$, $y_1$ and $y_2$, $|\Psi(x,y_1)|=|\Psi(x,y_2)|$, so that
$|\Psi(X_t)|$ is well-defined.  If we show that $\Psi(X_t,Y_t)$ is
with high probability close to $\E[\Psi(X_t,Y_t)]$, which in turn is
far from $0$, it will follow that $|\Psi(X_t)|$ is with high
probability confined to a small interval far from $0$, making it
statistically distinguishable from $|\Psi(X)|$ in stationarity,
implying that the Markov chain has not mixed by time $t$.

In the following sections we carry out these ideas to obtain lower bounds on
the mixing time that match (to within constants) the previously
obtained upper bounds.  Specifically, we show
\begin{theorem}
\label{thm:rudvalis}
For any fixed $\eps>0$, after
$\frac{1-o(1)}{8\pi^2} n^3\log n$ shuffles of the Rudvalis shuffle,
$\frac{1-o(1)}{2\pi^2} n^3\log n$ shuffles of the shift-or-swap shuffle, or
$\frac{1-o(1)}{\pi^2} n^3\log n$ shuffles for the symmetrized Rudvalis shuffle,
the distribution of the state of the deck has variation distance $\geq 1-\eps$
from uniformity.
\end{theorem}

\section{Lifting the Shuffles}
\label{sec:lift}

When the top card is placed at the bottom of the deck, the position of
any given card is cycically shifted left, so we will call this move
``shift-left'', and similarly ``shift-right'' is the move which places
the bottom card on top of the deck.  The move which exchanges the
top and bottom cards will be called ``swap'', and the move which does
nothing will be called ``hold''.  Thus the moves of the Rudvalis
Markov chain are ``shift-left'' and ``swap \& shift-left'', while the
moves of the variation considered by \dsc are ``shift-left'' and
``swap'', and the moves of the symmetrized version are ``shift-left'',
``shift-right'', ``swap'', and ``hold''.

The state $X_t$ of the Markov chain is the permutation giving the
order of the cards at time $t$.  Let $\diamondsuit$ (the diamond-suit
symbol) denote a particular card of interest, and let
$X_t(\diamondsuit)$ denote the location of card $\diamondsuit$ within
the deck at time $t$, where the positions are numbered from $1$ to $n$
starting from the top of the deck.  When the Markov chain does a
shift-left or shift-right, while the position of a card $\diamondsuit$
will change, all the cards get moved together, so it does not have
such a large randomizing effect on the permutation.  We will track the
position of card $\diamondsuit$, but we should also track the amount
of shifting.  So when we lift the Markov chain to $(X_t,Y_t)$, the
lifted Markov chain will also keeps track of
$$ Y_t = \text{\# shift-left's} - \text{\# shift-right's} \bmod n.$$
For the Rudvalis shuffle $Y_t=t\bmod n$ deterministically, whereas for
the other two variations, $Y_t$ will be a random number between $1$
and $n$ which approaches uniformity in $O(n^2)$ time.

Recall that we need an eigenvector $\Psi$ of the lifted chain $(X_t,Y_t)$
such that $|\Psi(X_t,Y_t)|$ is a function of $X_t$ alone.  For a given card
$\diamondsuit$, let
  $$\Psi_\diamondsuit(X_t,Y_t) = v(X_t(\diamondsuit)) \exp(Z_t(\diamondsuit) 2\pi i /n),$$
where $$Z_t(\diamondsuit) =
X_t(\diamondsuit)-X_0(\diamondsuit)+Y_t \bmod n,$$ and $v()$ is a
function, to be determined later, which makes $\Psi_\diamondsuit$ an
eigenvector.  Initially $Z_0(\diamondsuit)=0$, and the only time that
the $Z_t(\diamondsuit)$ changes is when the card $\diamondsuit$ gets
transposed.
The dynamics of $(X_t(\diamondsuit),Z_t(\diamondsuit))$ (mod $n$) are summarized by
$$ (X_{t+1}(\diamondsuit),Z_{t+1}(\diamondsuit)) = \begin{cases}
   (X_t(\diamondsuit),Z_t(\diamondsuit)) & \text{if move was ``hold''} \\
   (X_t(\diamondsuit)-1,Z_t(\diamondsuit)) & \text{if move was ``shift-left''} \\
   (X_t(\diamondsuit)+1,Z_t(\diamondsuit)) & \text{if move was ``shift-right''} \\
   (X_t(\diamondsuit)-1,Z_t(\diamondsuit)-1) & \text{if move was ``swap'' and $X_t(\diamondsuit)=1$} \\
   (X_t(\diamondsuit)+1,Z_t(\diamondsuit)+1) & \text{if move was ``swap'' and $X_t(\diamondsuit)=n$} \\
   (X_t(\diamondsuit),Z_t(\diamondsuit)) & \text{if move was ``swap'' and $\diamondsuit$ elsewhere} \\
 \end{cases}
$$
We define
 $$\Psi(X_t,Y_t) = \sum_{\diamondsuit=1}^{n} \Psi_\diamondsuit(X_t,Y_t).$$
If we increment $y$ while holding $x$ fixed, then $\Psi(x,y)$
gets multiplied by the phase factor $\exp(2\pi i/n)$, so we have an
eigenvector satisfying our requirement that $|\Psi(X_t,Y_t)|$ be a
function of $X_t$ alone.

\section{The Lower Bound Lemma}

The lower bounding lemma that we shall use is similar to Lemma~4 of
\citep{wilson:mixing-times}, but with the modifications described in
the introduction.  \cite{saloff-coste:lower-bound} also
gives a generalization of Lemma~4 from \citep{wilson:mixing-times} that
may be used when the eigenvalues are complex, but the extension below
seems to be better suited for the shuffles considered here.
\begin{lemma}\label{lem:anticonverge}
Suppose that a Markov chain $X_t$ has a lifting $(X_t,Y_t)$, and that
$\Psi$ is an eigenfunction of the lifted Markov chain:
$\E[\Psi(X_{t+1},Y_{t+1})\mid (X_t,Y_t)] = \lambda \Psi(X_t,Y_t)$.
Suppose that $|\Psi(x,y)|$ is a function of $x$ alone, $|\lambda|<1$,
$\Re(\lambda)\geq 1/2$, and that we have an upper bound $R$ on
$\E[|\Psi(X_{t+1},Y_{t+1})-\Psi(X_t,Y_t)|^2\mid (X_t,Y_t)]$.
Let $\gamma=1-\Re(\lambda)$.
Then when the number of Markov chain steps $t$ is bounded by
$$ t \leq \frac{\log \px + \frac{1}{2} \log \frac{\gap \varepsilon }{4 R}
  }{- \log(1-\gap)},$$
the variation distance of $X_t$ (the state of the original Markov chain)
from stationarity is at least $1-\varepsilon$.
\end{lemma}

The proof of this modified lemma is similar to the proof of Lemma~4 in
\citep{wilson:mixing-times}, but for the reader's convenience we give
the modified proof.  In the following sections we determine the
functions $v()$ for the Markov chains which give the requisite
eigenfunction $\Psi$, and then use \lref{anticonverge} to obtain the
mixing time bounds stated in \tref{rudvalis}.

\begin{proof}[Proof of \lref{anticonverge}]
Let $\Psi_t = \Psi(X_t,Y_t)$, and $\dps=\Psi_{t+1}-\Psi_t$.
By induction $$\E[\Psi_t\mid (X_0,Y_0)] = \Psi_0 \lambda^t.$$
By our assumptions on $\lambda$, in equilibrium $\E[\Psi]=0$.

We have $\E[\dps\mid (X_t,Y_t)]=(\lambda-1)\Psi_t$ and
\begin{align*}
\Psi_{t+1}\Psi_{t+1}^*  &= \Psi_t\Psi_t^* + \Psi_t\dps^* + \Psi_t^*\dps + |\dps|^2\\
\E[\Psi_{t+1}\Psi_{t+1}^*\mid (X_t,Y_t)]&= \Psi_t \Psi_t^*[1+(\lambda-1)^*+(\lambda-1)] + \E[|\dps|^2\mid X_t] \\
 &\leq 
\Psi_t \Psi_t^*[2\Re(\lambda)-1] + R \\
\intertext{and so by induction,}
\E[\Psi_t\Psi_t^*]  &\leq \Psi_0\Psi_0^* [2\Re(\lambda)-1]^t + \frac{R}{2-2\Re(\lambda)},
\intertext{then subtracting $\E[\Psi_t]\E[\Psi_t]^*$,}
\Var[\Psi_t]   &\leq \Psi_0\Psi_0^* \left[[2\Re(\lambda)-1]^t - (\lambda\lambda^*)^t\right]
                  + \frac{R}{2-2\Re(\lambda)} .\\
\intertext{Since $(1-\lambda)(1-\lambda^*)\geq 0$, we have $\lambda\lambda^* \geq 2\Re(\lambda)-1$, and by assumption $2\Re(\lambda)-1\geq 0$.  Hence $(\lambda\lambda^*)^t \geq [2\Re(\lambda)-1]^t$, and we have for each $t$}
\Var[\Psi_t]   &\leq \frac{R}{2-2\Re(\lambda)} = \frac{R}{2\gamma}.
\end{align*}

From Chebychev's inequality, $$\Pr\left[|\Psi_t-\E[\Psi_t]| \geq
  \sqrt{R/(2\gap\varepsilon)}\right] \leq \varepsilon.$$
As $\E[\Psi_\infty]=0$, if $\E[\Psi_t] \geq \sqrt{4 R/(\gap\varepsilon)}$,
then the probability that $|\Psi_t|$ deviates below
$\sqrt{R/(\gap\varepsilon)}$ is at most $\varepsilon/2$, and the
probability that $|\Psi|$ in stationarity deviates above this threshold
is at most $\varepsilon/2$, so the variation distance between the
distribution at time $t$ and stationarity must be at least
$1-\varepsilon$.
If we take the initial state to be the one maximizing $\Psi_0$,
then $$\E[|\Psi_t|] = |\px| |\lambda|^t \geq |\px| (\Re(\lambda))^t = |\px| (1-\gamma)^t \geq \sqrt{4 R/(\gap\varepsilon)}$$
when
\begin{align*}
 t &\leq \frac{\log \left[\px\div\sqrt{\frac{4 R}{\gap\varepsilon}}\right]}{-\log(1-\gap)}.\qedhere
\end{align*}
\end{proof}

\section{The Rudvalis Shuffle}
\label{sec:rudvalis}

The first shuffle we consider is the original Rudvalis Markov chain.
It will be instructive to consider a slight generalization, where the
swap \& shift-left move takes place with probability $p$, and the
shift-left move takes place with probability $1-p$.  We shall assume
that $0<p<1$ and that $p$ is independent of $n$.  The particular
values of $p$ that we are interested in are $p=1/2$ (for the original
Rudvalis chain) and $p=1/3$.

We need to find an eigenvector for the random walk that a single card
takes under this shuffle.  We remark that this random walk is similar
in nature to (but distinct from) a class of random walks, known as
daisy chains, for which \cite{wilmer:thesis} obtained eigenvalues and
eigenvectors.  From other work of \cite{wilmer:necklace}, it readily
follows that the position of a single card takes order $n^3$ steps to
randomize, and that the precise asymptotic distance from stationarity
of the card's position after $c n^3$ shuffles is given by an explicit
expression involving theta functions.

\begin{lemma}
The random walk followed by a card $\diamondsuit$ under the lifted
Rudvalis shuffle has an eigenvector of the form
$$ \Psi_\diamondsuit(x,z) = v(x) e^{2\pi i z/n} $$
where $v(x)$ is the $x^{\text{th}}$ number in the list
$$\lambda^{n-2},\ldots,\lambda^2,\lambda,1,\chi\ ,$$
the eigenvalue is
$$\lambda = 1 - \frac{p}{1-p} \frac{4\pi^2}{n^3} + O(1/n^4),$$
and
$$\chi=1+\frac{p}{1-p} \frac{2\pi i}{n} + O(1/n^2).$$
\end{lemma}

\begin{proof}
Let $w = \exp(2\pi i/n)$.
If at time $t$ card $\diamondsuit$ is in any location between $2$ and $n-1$, then
$$\Psi_\diamondsuit(X_{t+1},Y_{t+1}) = \lambda \Psi_\diamondsuit(X_t,Y_t)$$
deterministically.  To ensure that 
$$\E[\Psi_\diamondsuit(X_{t+1},Y_{t+1})\mid (X_t,Y_t)] = \lambda \Psi_\diamondsuit(X_t,Y_t)$$
when $X_t(\diamondsuit) = 1$, we require
\begin{align*}
p w^{-1} + (1-p) \chi &= \lambda^{n-1} \\
          \chi &= \frac{\lambda^{n-1} - p w^{-1}}{1-p},
\end{align*}
and for when $X_t(\diamondsuit) = n$ we need
\begin{align*}
 p w \chi + (1-p) &= \lambda \chi \\
           \chi &= \frac{1-p}{\lambda-p w}.
\end{align*}
Given these two equations, $\Psi_\diamondsuit$ will
be an eigenvector with eigenvalue $\lambda$.
Thus,
$$
f(\lambda) = \lambda^n - p w \lambda^{n-1} - p w^{-1} \lambda - 1 + 2 p = 0.
$$

To identify a root of this polynomial,
we use Newton's method:
$z_{k+1} = z_k - f(z_k)/f'(z_k)$, starting with $z_0=1$.
By Taylor's theorem,
$$|f(z_{k+1})| \leq \frac{1}{2} \max_{0\leq u \leq 1} |f''(u z_k + (1-u) z_{k+1})| \times \left|\frac{f(z_k)}{f'(z_k)}\right|^2.$$
If $|z-1|\leq 1/n^2$, then
$f'(z) = (1-p) n + O(1)$ and $f''(z) = (1-p) n^2 + O(n)$.
Consequently, if $|z_k-1|\leq 1/n^2$ and $|z_{k+1}-1|\leq 1/n^2$,
then $$|f(z_{k+1})| \leq \frac{1+O(1/n)}{2} \frac{1}{1-p} |f(z_k)|^2.$$
Since $f(z_0)=p(2-w-w^{-1}) = p 4\pi^2/n^2 + O(1/n^4)$, for large enough
$n$ we have by induction that $|f(z_k)| \leq (1-p) (p/(1-p) 4\pi^2/n^2)^{2^k}$,
$|z_{k+1}-z_k| \leq (p/(1-p) 4\pi^2/n^2)^{2^k}/(n+O(1))$, and $|z_{k+1}-z_0|\leq O(1/n^3)$.
Thus, for large enough $n$, the sequence $z_0,z_1,z_2,\ldots$ converges to
a point $\lambda$, which by continuity, must be a zero of $f$.  Since
$z_1 = 1-p/(1-p) 4\pi^2/n^3 + O(1/n^4)$ and $|\lambda-z_1|=O(1/n^5)$,
we conclude that the polynomial $f$ has a root at
\begin{align*}
\lambda &= 1 - \frac{p}{1-p} \frac{4\pi^2}{n^3} + O(1/n^4) .\qedhere
\end{align*}
\end{proof}

It is easy to check that $\px=n+O(1/n)$.
Next we evaluate $R$ for this eigenvector.
$$
\frac{\Psi_\diamondsuit(X_{t+1},Y_{t+1}) - \Psi_\diamondsuit(X_t,Y_t)}
{w^{Z_t(\diamondsuit)}} =
\begin{cases}
(\lambda-1)\lambda^{X_t(\diamondsuit)} = O(1/n^3)
 & \text{if $2\leq X_t(\diamondsuit)\leq n-1$} \\
 \chi - \lambda^{n-2} = O(1/n)
 & \text{if $X_t(\diamondsuit)=1$, shift-left} \\
 w^{-1} - \lambda^{n-2} = O(1/n)
 & \text{if $X_t(\diamondsuit)=1$, swap \& shift-left} \\
 1-\chi = O(1/n)
 & \text{if $X_t(\diamondsuit)=n$, shift-left} \\
 w\chi-\chi = O(1/n)
 & \text{if $X_t(\diamondsuit)=n$, swap \& shift-left} \\
 \end{cases}
$$
Adding up these contributions over the various cards $\diamondsuit$, we find
\begin{align*}
|\Psi(X_{t+1},Y_{t+1})-\Psi(X_t,Y_t)| &\leq O(1/n) \\
R = \E[|\Psi(X_{t+1},Y_{t+1})-\Psi(X_t,Y_t)|^2\mid (X_t,Y_t)] &\leq O(1/n^2).
\end{align*}
Plugging $\lambda$, $\px$, and $R$ into the \lref{anticonverge} gives, for fixed values of $\eps$,
 a mixing time lower bound of
$$ (1-o(1))\frac{1-p}{p} \frac{1}{8\pi^2} n^3 \log n.$$

\section{The Shift-or-Swap Shuffle}
\label{sec:shift-or-swap}

At this point there are two ways we can approach the shift-or-swap
shuffle.  We can either take a direct approach in the same manner as
in the previous section, or we can do a comparison with the Rudvalis
shuffle with $p=1/3$.

If we take the direct approach, then we let $v(x)$ denote the
$x^{\text{th}}$ element of the list
$$ (2\lambda-1)^{n-2},\ldots,(2\lambda-1)^2,2\lambda-1,1,\chi.$$
The constraints on $\chi$ are
$$ \chi = \frac{2\lambda}{1+w^{-1}}(2\lambda-1)^{n-2} $$
and
$$ \chi = \frac{1+w(2\lambda-1)^{n-2}}{2\lambda}. $$
As in section~\ref{sec:rudvalis}, we solve for $\lambda$ and find that
$\lambda = 1-(1+o(1))\pi^2/n^3$, compute $\px=\Theta(n)$ and $R=O(1/n^2)$,
and obtain the mixing time lower bound of $\frac{1-o(1)}{2\pi^2} n^3\log n$ shuffles.

Alternatively, we can couple the shift-or-swap shuffle with the
Rudvalis shuffle.  Whenever the shift-or-swap shuffle makes a shift,
the number of swap's since the previous shift will be odd with
probability $1/3$.  If it is odd, then this is equivalent to a
swap-\&-shift-left move, and if it is even, then it is equivalent to a
shift-left move.  This explains why we were interested in the case
$p=1/3$ in the previous section.  After $t$ steps, with high
probability $(1+o(1)) t/2$ shift moves occured, which means that the
state of the deck is what it would be after $(1+o(1)) t/2$ Rudvalis
shuffles (with $p=1/3$), possibly with an extra swap move.  The lower
bound for the shift-or-swap shuffle does not follow from the lower
bound itself for the Rudvalis shuffle, but it does follow from what
we showed about $|\Psi|$ for the Rudvalis shuffle.

\section{Symmetrized Version of the Rudvalis Shuffle}
\label{sec:symmetrized}

When analyzing the symmetrized version of the Rudvalis shuffle,
it will be convenient to have symmetric coordinates, so we
re-index the card locations to run from $-(n-1)/2$ up to $(n-1)/2$, and the
swaps occur at locations $-(n-1)/2$ and $(n-1)/2$.

\begin{lemma}
The random walk followed by a card $\diamondsuit$ under the lifted
symmetrized Rudvalis shuffle has an eigenvector of the form
$$ \Psi_\diamondsuit(x,z) = v(x) e^{2\pi i z/n} $$
where
$$v(x) = \frac{1+\delta}{2} e^{i\theta x} + \frac{1-\delta}{2} e^{-i\theta x}
       = \cos(\theta x)+i\delta \sin(\theta x),$$
$$\theta = (1+o(1)) \sqrt{2}\pi n^{-3/2},$$
$$ \delta = (1+o(1))\frac{1}{\sqrt{2} n^{1/2}},$$
and both $\delta$ and $\theta$ are real.
The eigenvalue is
$$  \lambda = \frac{1+\cos\theta}{2} = 1-\frac{\pi^2+o(1)}{2} n^{-3}. $$
\end{lemma}

\begin{proof}
When $x\neq\pm(n-1)/2$, we can readily compute the eigenvalue $\lambda$
to be
\begin{align*}
 \lambda &= \frac{\frac14 v(x+1)+\frac12 v(x)+\frac14 v(x-1)}{v(x)} \\
     &= \frac12 \frac{\cos(\theta x)\cos\theta + i\delta\sin(\theta x)\cos\theta}{\cos(\theta x)+i\delta \sin(\theta x)} +\frac12\\
 &= \frac{1+\cos\theta}{2}.
\end{align*}
In order for our guessed eigenvector to be correct,
there is also a constraint at $x=(n-1)/2$:
\begin{align*}
\frac{v(\frac{n-1}{2})}{4} + \frac{v(\frac{n-3}{2})}{4} + (1+w)\frac{v(-\frac{n-1}{2})}{4}
 &= \lambda = \frac{v(\frac{n-1}{2})}{2} + \frac{v(\frac{n-3}{2})}{4} + \frac{v(\frac{n+1}{2})}{4} \\
(1+w) v(-(n-1)/2) &= v((n-1)/2) + v((n+1)/2) \\
(1+w) (1+\delta) e^{-i\theta (n-1)/2} + (1+w)(1-\delta) e^{i\theta (n-1)/2}
      &=
\begin{aligned}[t] 
& (1+\delta) e^{i\theta (n-1)/2} + (1-\delta) e^{-i\theta (n-1)/2} +\\&
  (1+\delta) e^{i\theta (n+1)/2} + (1-\delta) e^{-i\theta (n+1)/2}
\end{aligned} \\
(w+2\delta+w\delta) e^{-i\theta (n-1)/2} + (w-2\delta-w\delta) e^{i\theta (n-1)/2}
      &= 
(1+\delta) e^{i\theta (n+1)/2} + (1-\delta) e^{-i\theta (n+1)/2} \\
\textstyle
w\cos\frac{\theta (n-1)}{2} - \cos\frac{\theta (n+1)}{2} 
      &= \textstyle
 \delta \left[ (2+w) i \sin\frac{\theta (n-1)}{2} + i \sin\frac{\theta (n+1)}{2}\right].
\end{align*}
The corresponding constraint at $x=-(n-1)/2$ is obtained by replacing
$w$ with $1/w$ and replacing $\theta$ with $-\theta$.  Since these
substitutions give the complex-conjugate of the above equation, the
constraints at $x=\pm(n-1)/2$ are equivalent.

Equating the real parts of this equation gives
\begin{align*}
 \delta &= \frac{\cos\frac{2\pi}{n} \cos\frac{\theta (n-1)}{2} - \cos\frac{\theta (n+1)}{2}}
              {-\sin\frac{2\pi}{n} \sin\frac{\theta (n-1)}{2}}, \\
\intertext{and equating the imaginary parts gives}
 \delta &= \frac{\sin\frac{2\pi}{n} \cos\frac{\theta (n-1)}{2}}{2 \sin\frac{\theta (n-1)}{2} + \cos\frac{2\pi}{n} \sin\frac{\theta (n-1)}{2} + \sin\frac{\theta (n+1)}{2}} .
\end{align*}
Cross-multiplying and performing trigonometric simplifications gives
$$
 -\sin^2(\textstyle\frac{2\pi}{n}) \displaystyle\frac{\sin(\theta (n-1))}{2} =
\begin{aligned}[t]
 & \cos^2(\textstyle\frac{2\pi}{n}) \displaystyle\frac{\sin(\theta (n-1))}{2} - \displaystyle\frac{\sin(\theta (n+1))}{2} + \cos(\textstyle\frac{2\pi}{n}) \sin\theta \\ &+ \cos(\textstyle\frac{2\pi}{n})\sin(\theta(n-1))-2\sin\textstyle\frac{\theta (n-1)}{2} \cos\textstyle\frac{\theta (n+1)}{2}
\end{aligned}
$$
so that
\begin{equation}
\label{eq:theta}
  \left(\frac12+\cos\frac{2\pi}{n}\right)\sin(\theta (n-1)) - \frac12 \sin(\theta (n+1)) -\sin(\theta n) +\left(1+\cos\frac{2\pi}{n}\right)\sin(\theta) = 0.
\end{equation}
Equation~\eqref{eq:theta} is exact, but to estimate a solution, we perform
a series expansion in $\theta$
\begin{align*}
 0 &= \begin{aligned}[t]
      &\left[\left(\frac12+\cos\frac{2\pi}{n}\right) (n-1)-(n+1)/2-n+1+\cos\frac{2\pi}{n}\right] \theta \\ &-
      \left[\left(\frac12+\cos\frac{2\pi}{n}\right) (n-1)^3-(n+1)^3/2-n^3 +1+\cos\frac{2\pi}{n}\right]\frac{\theta^3}{6} + O(n^4\theta^5)
      \end{aligned} \\
  0 &= -\left[\frac{2\pi^2}{n}+O(1/n^2)\right] \theta + \left[6n^2+O(n)\right]\frac{\theta^3}{6} + O(n^4\theta^5).
\end{align*}
While $\theta=0$ is a solution, our expression for $\delta$ has a singularity at $\theta=0$, so we seek a different solution.  Ignoring the error terms suggests $\theta \doteq \sqrt{2}\pi n^{-3/2}$.
Since the function in \eqref{eq:theta} is real-valued, we can appeal to the intermediate value theorem to show that there is in fact a root at 
$$\theta\doteq \sqrt{2}\pi n^{-3/2}.$$
For this value of $\theta$ we have
$$  \lambda = \frac{1+\cos\theta}{2} \doteq 1-\frac{\pi^2}{2} n^{-3}, $$
and (using the second equation for $\delta$)
\begin{align*}
 \delta &\doteq \frac{2\pi/n}{2 n \theta/2 + n\theta/2 + n\theta/2} \doteq \frac{1}{\sqrt{2} n^{1/2}}.\qedhere
\end{align*}
\end{proof}

Again $\px=(1+o(1))n$.
Next we estimate $R$.  If there is a shift-left, then provided card~$\diamondsuit$ is not in position $-(n-1)/2$, we have
\begin{align*}
 \Delta\Psi_\diamondsuit
 &= (\cos\theta-1)[\cos(\theta x)+i\delta\sin(\theta x)]e^{2\pi i z/n} + \sin\theta[\sin(\theta x)-i\delta\cos(\theta x)]e^{2\pi i z/n} \\
 &= O(\theta^2) + O(\theta^2 x) + O(\theta\delta) = O(n^{-3})+O(n^{-2})+O(n^{-2}) = O(n^{-2}).
\end{align*}
If card $\diamondsuit$ is in position $-(n-1)/2$, then
$$
 \Delta\Psi_\diamondsuit =
  2 i \delta \sin(\theta(n-1)/2) e^{2\pi i z/n} = O(n^{-1}).
$$
Adding up these contributions over the different cards, we find
$ \Delta\Psi = O(n^{-1})$.
Likewise $\Delta\Psi=O(n^{-1})$ if the move was a shift-right.  For
transposes, $\Delta\Psi_\diamondsuit$ is nonzero for only two cards,
and for these it is $O(n^{-1})$.  Thus in all cases we have
$|\Delta\Psi|^2 \leq O(n^{-2})$, and so $R \leq O(n^{-2})$.
Plugging our values of $\lambda$, $\px$, and $R$ into \lref{anticonverge}, we
obtain, for fixed values of $\eps$, a lower bound on the mixing time of
$ \frac{1-o(1)}{\pi^2} n^3 \log n$ shuffles.

\section{Remarks}

We have seen how to extend the lower bound technique used in
\citep{wilson:mixing-times} to % non-reversible shuffles and
shuffles
that are much slower than what the position of a single card would
indicate.  Interestingly, for the shift-or-swap and
symmetrized-Rudvalis shuffles, the spectral gap for the lifted shuffle
is smaller than the spectral gap of the shuffle itself, so it is
curious that we obtained a lower bound for these shuffles by
considering their lifted versions.

In an early draft, we lower bounded the mixing time of the original
Rudvalis shuffle without considering its lifted version, and this
earlier approach might be considered simpler.  But it is not clear how
to lower bound the symmetrized Rudvalis shuffle without lifting it, and
our current approach has the advantage that the analyses for all
three shuffles treated here are similar.  The original Rudvalis
shuffle and its lifting are isomorphic, and the earlier analysis
is effectively a special case of the present analysis
where the lifting is not explicit.

We suspect that the constants in the lower bounds of
\tref{rudvalis} are tight.  Rudvalis asked if his shuffle was
the slowest shuffle evenly supported on two generators; the lower
bounds given here suggest that the shift-or-swap shuffle (for odd $n$)
is slower by a factor of $4$.

\section*{Acknowledgements}

The author thanks Laurent Saloff-Coste and Persi Diaconis for calling
attention to these shuffles, and Laurent Saloff-Coste for comments on
an earlier draft.

\bibliography{rudvalis}

\begin{thebibliography}{10}
\expandafter\ifx\csname natexlab\endcsname\relax\def\natexlab#1{#1}\fi

\bibitem[Aldous and Fill(2005)]{aldous-fill:book}
David~J. Aldous and James~A. Fill.
\newblock {\em {Reversible {M}arkov Chains and Random Walks on Graphs}}.
\newblock Book in preparation, {\tt
  http://www.stat.berkeley.edu/\char126aldous/book.html}, 2005.

\bibitem[Diaconis(1988)]{diaconis:monograph}
Persi Diaconis.
\newblock {\em Group Representations in Probability and Statistics}.
\newblock Institute of Mathematical Statistics, 1988.

\bibitem[Diaconis(1996)]{diaconis:cutoff}
Persi Diaconis.
\newblock The cutoff phenomenon in finite {Markov} chains.
\newblock {\em Proceedings of the National Academy of Sciences, USA},
  93:\penalty0 1659--1664, 1996.

\bibitem[Diaconis and Saloff-Coste(1993)]{diaconis-saloff-coste:comparison}
Persi Diaconis and Laurent Saloff-Coste.
\newblock Comparison techniques for random walk on finite groups.
\newblock {\em The Annals of Probability}, 21\penalty0 (4):\penalty0
  2131--2156, 1993.

\bibitem[Diaconis and Saloff-Coste(1995)]{diaconis-saloff-coste:survey}
Persi Diaconis and Laurent Saloff-Coste.
\newblock Random walks on finite groups: a survey of analytic techniques.
\newblock In {\em Probability measures on groups and related structures, XI
  (Oberwolfach, 1994)}, pages 44--75. World Sci. Publishing, 1995.

\bibitem[Hildebrand(1990)]{hildebrand:thesis}
Martin~V. Hildebrand.
\newblock {\em Rates of Convergence of Some Random Processes on Finite Groups}.
\newblock PhD thesis, Harvard University, 1990.

\bibitem[Saloff-Coste(2002)]{saloff-coste:lower-bound}
Laurent Saloff-Coste.
\newblock Lower bound in total variation for finite {Markov} chains: {Wilson's}
  lemma, 2002.
\newblock Manuscript.

\bibitem[Wilmer(1999)]{wilmer:thesis}
Elizabeth~L. Wilmer.
\newblock {\em Exact Rates of Convergence for Some Simple Non-Reversible
  {Markov} Chains}.
\newblock PhD thesis, Harvard University, 1999.

\bibitem[Wilmer(2002)]{wilmer:necklace}
Elizabeth~L. Wilmer.
\newblock A local limit theorem for a family of non-reversible {Markov} chains,
  2002.
\newblock arXiv:math.PR/0205189.

\bibitem[Wilson(2001)]{wilson:mixing-times}
David~B. Wilson.
\newblock Mixing times of lozenge tiling and card shuffling {Markov} chains,
  2001.
\newblock To appear in \textit{The Annals of Applied Probability}.
  arXiv:math.PR/0102193.

\end{thebibliography}
\bibliographystyle{plainnat}

\end{document}